\documentclass{singlecol-new}

\usepackage{natbib,stfloats}
\usepackage{mathrsfs,amsmath,upgreek}

\newcommand{\be}{\begin{eqnarray}}
\newcommand{\ee}{\end{eqnarray}}

\theoremstyle{TH}{
\newtheorem{lemma}{Lemma}
\newtheorem{theorem}[lemma]{Theorem}

\newtheorem{proposition}[lemma]{Proposition}

}
\theoremstyle{THrm}{
\newtheorem{definition}{Definition}

\newtheorem{remark}{Remark}

}

\theoremstyle{THhit}{

}

\makeatletter

\def\theequation{\arabic{equation}}

\def\tc{\textcolor{red}}

 \makeatother

\begin{document}%

\thispagestyle{plain}

\setcounter{page}{1}

\LRH{}

\RRH{Clustering and Hitting Times of Threshold Exceedances and  Applications}

\VOL{x}

\ISSUE{x}

\PUBYEAR{2015}




\title{Clustering and Hitting Times of Threshold Exceedances and  Applications}

\authorA{Natalia Markovich}

\affA{V. A. Trapeznikov Institute of Control Sciences
\\
of Russian Academy of Sciences,
\\
Profsoyuznaya str. 65,
Moscow 117997, Russia
\\
Email: nat.markovich@gmail.com, markovic@ipu.rssi.ru}


\begin{abstract}
We investigate exceedances of the process over a sufficiently high threshold. The exceedances determine the risk of hazardous events like climate catastrophes, huge insurance claims, the loss and delay in telecommunication networks.
 Due to dependence such exceedances tend to occur in clusters. The cluster structure of social networks is caused by dependence (social relationships and
interests) between nodes and possibly heavy-tailed distributions of the node
degrees. A minimal time to reach a large node determines the first hitting time.
We derive an asymptotically equivalent distribution and a limit expectation of the first
hitting time to exceed the threshold $u_n$ as the sample size $n$ tends to infinity. The results can be extended to the second and, generally, to the $k$th ($k> 2$) hitting times. Applications in large-scale networks such as social, telecommunication and recommender systems are discussed.
\end{abstract}

\KEYWORD{first hitting time; rare events; exceedance over threshold; cluster of exceedances; extremal index; application.}

\REF{to this paper should be made as follows:  (xxxx) `Clustering and Hitting Times of Threshold Exceedances and  Application',
{\it International Journal of Data Analysis Techniques and Strategies}, Vol.~x, No.~x, pp.xxx--xxx.}

\begin{bio}
Natalia Markovich received the M.S. degree in Mathematics (with distinction) from the  Applied Mathematics and Cybernetics Department of Tver University,  Russia,   Ph.D. and Doctor Sciences degree  in Statistics from the Institute of Control Sciences of Russian Academy of Sciences in Moscow. She is a Main Scientist in the Institute of Control Sciences, and  her research interests concern extreme value theory,
nonparametric statistics,  statistical analysis of measurements in telecommunication systems, including the Internet.  She has published over 60 papers and the book "Nonparametric analysis of univariate heavy-tailed data. Research and Practice", Wiley, 2007.
\end{bio}

\maketitle

\section{Introduction}
\label{sec-introduction}
Let $\{X_n\}_{n\ge 1}$ be a stationary sequence with marginal distribution function $F(x)$ and $M_n=\max\{X_1,...,X_n\}$. We investigate rare events, namely,  exceedances of the sequence over a sufficiently high threshold $u$. Due to dependence such exceedances tend to occur in clusters. Such clusters of rare events and the asymptotic distributions of the cluster and inter-cluster sizes have been widely studied due to numerous applications, see \cite{ancona00}, \cite{beirlant}, \cite{ferro}, \cite{markovich14}, \cite{markovich16a}, \cite{robert09}, \cite{robert13}, \cite{roberts}, \cite{robinson} among others. There are three approaches in the cluster size study, namely, the blocks method, the runs method and the inter-exceedance times method. The first two methods define the cluster as a block of data with at least one exceedance over the threshold or the clusters are blocks of data with some number of exceedances which are separated  by at least a fixed number of observations running under the threshold, respectively, \cite{SmithWeissman}, \cite{WeissmanNovak}. Following the inter-exceedance times approach proposed in \cite{ferro} we define the cluster as a conglomerate 
of consecutive exceedances over the threshold between two consecutive non-exceedances. Our main objective is to study the distribution of the first hitting time to exceed the threshold $u$.
\\
 Let us consider the inter-cluster size
 \begin{eqnarray}\label{18}T_1(u)&=&\min\{j\ge 1: M_{1,j}\le u,X_{j+1}>u|X_{1}>u\},\end{eqnarray} i.e. the number of inter-arrivals of observations running under the threshold between two consecutive exceedances,
where $M_{1,j}=\max\{X_{2},...,X_{j}\}$, $M_{1,1}=-\infty$.
Let \[T^*(u)=\min\{j+1\ge 1: M_{j}\le u,X_{j+1}>u\}\] be the first hitting time corresponding to the threshold $u$. We get
\begin{eqnarray}\label{0}P\{T^*(u)=j+1\}&=& P\{M_j\le u, X_{j+1}>u\},\end{eqnarray} $j=0,1,2,...$, $M_{0}=-\infty$.
\\
Let $T^*_T(u)$ be the first hitting time in the  time interval $[0,T]$. 
Let $\{Y_n\}_{n\ge 1}$ be a stationary sequence of inter-arrival times between consecutive observations of the $\{X_n\}$ and $S_j=\sum_{i=1}^{j-1} Y_i$ denotes the time interval between arrivals of $X_1$ and $X_j$.
Then we have
\begin{eqnarray}\label{10a}P\{T^*_T(u)=j+1\}&=& P\{M_j\le u, X_{j+1}>u, S_{j+1}\le T\}.\end{eqnarray}
Similarly, we determine the probability of the $k$ consecutive hitting times $T^{**}(u)$ by
\begin{eqnarray*}\label{12}P\{T^{**}(u)=k\}&=& P\{M_{i_1}\le u, X_{i_1+1}>u, M_{i_1+1,i_2}\le u, X_{i_2+1}>u,\nonumber
\\
&...,& M_{i_{k-1}+1,i_k}\le u, X_{i_k+1}>u\},\end{eqnarray*}
$i_j=0,1,2,...$; $j=1,2,...,k$.
\\
The necessity to evaluate the distribution, quantiles and the mean of the  first hitting time is arising in many applications. In social networks it is important to compare sampling strategies (\cite{avrachenkov12, avrachenkov15, lee}) like  random walks, Metropolis-Hastings Markov chains, Page Ranks and others with regard to how quickly they allow to reach a node with a large degree, that is the number of links with other nodes. In \cite{markovich15} it is proposed to compare sampling techniques by the mean first hitting time that is illustrated on the real data of social networks. It is important to investigate the first hitting time of significant nodes since it allows us to disseminate advertisements or to collect opinions more effectively within clusters surrounding such nodes. It can be helpful also for recommender systems with collaborative filtering, in which the system recommends to a user some item or product  that has been rated by previous users, 
\cite{lu}. A similar problem occurs in telecommunication peer-to-peer networks, namely to find a node with a large number of peers, \cite{dan, markovich13}. Our concept is also relevant in other areas of operations research and inventory control. For instance, the first hitting time can be important to analyze the customer churn that has a huge impact on companies, \cite{Mahajan}. Considering the customer impatience in multi-server queues \cite{Choudhury} and the customer waiting time in the queue \cite{Zhao}, the first hitting time indicates the moment when the waiting time  exceeds a threshold and hence, the customer may leave the queue.  Following \cite{Lattila} the forecasting of exceedances of industrial production can be a driving factor for the development of sea ports.
The $k$th hitting time is important in Internet to find the $k$-top web sites that are significant with regard to some topic.
\\
The measure of the dependence between the rare events is expressed by the extremal index.
The notion of the extremal index is determined in \cite{leadbetter83}, p.53.
\begin{definition} The stationary sequence  $\{X_n\}_{n\ge 1}$ is said to have extremal index
$\theta\in[0,1]$ if
for each $0<\tau <\infty$ there is a sequence of real numbers $u_n=u_n(\tau)$ such that
\begin{equation}\label{1}\lim_{n\to\infty}n(1-F(u_n))=\tau \qquad\mbox{and}\end{equation}
\begin{equation}\label{2}\lim_{n\to\infty}P\{M_n\le u_n\}=e^{-\tau\theta}\end{equation}
hold.
\end{definition}
The  extremal index $\theta$ of  $\{X_n\}$ relates to the first hitting time  $T^*(u_n)$, \cite{roberts}.
Really, since
 $u_n$ is selected according to (\ref{1}) it follows that $P\{X_n > u_n\}$ is asymptotically equivalent to $1/n$. Notice, that $P\{M_k\le u_n\}=P\{T^*(u_n)>k\}$. Hence, substituting $\tau$ by (\ref{1})  we get from (\ref{2}) \[P\{T^*(u_n)/n>k/n\}\sim e^{-\theta k P\{X_n > u_n\}}\sim e^{-\theta k/n},\] 
 \[\lim_{n\to\infty}P(T^*(u_n)/n > x) = e^{-\theta x}\] for positive $x$. It follows
 \begin{equation}\label{5bc}\lim_{n\to\infty}E(T^*(u_n)/n)=1/\theta.\end{equation}
 This implies,  the smaller $\theta$, the longer it takes to reach an observation with a large value. The result is then interesting for processes which have $\theta=0$, e.g., for the  Metropolis Markov chain \cite{roberts} and the Lindley process with subexponential step distribution \cite{asmussen1998}. A mixture of i.i.d. non-ergodic sequences with $\theta=0$ is given in Theorem 4 by \cite{doukhan}.
\\
Using achievements regarding the limit geometric-like distribution of $T_1(x_{\rho_n})$ derived in (Theorem 2, \cite{markovich14}, \cite{markovich16a}), where the $(1-\rho_n)$th quantile $x_{\rho_n}$ of $\{X_n\}$ is taken as $u_n$, we derive in Section \ref{Sec2} a limit distribution of the first hitting time and its expectation that specifies (\ref{5bc}). The achievements are similarly extended to the second 
hitting time, Section \ref{Sec3}.
\\
Theorem \ref{T2} (that is Theorem 2 in \cite{markovich14})
is based on the mixing condition proposed in \cite{ferro}
\begin{equation}\label{44}\alpha_{n,q}(u)=\max_{1\le k\le n-q}\sup|P(B|A)-P(B)|=o(1),\qquad n\to\infty,
\end{equation}where for real  $u$ and integers $1\le k\le l$,
$\mathcal{F}_{k,l}(u)$
is the $\sigma$-field generated by the events $\{X_i>u\}$, $k\le
i\le l$ and the supremum is taken over all $A\in
\mathcal{F}_{1,k}(u)$ with $P(A)> 0$ and $B\in
\mathcal{F}_{k+q,n}(u)$ and $k$, $q$ are positive integers.
\\
To formulate the theorem we need  the following partition of the interval $[1,j]$
\begin{equation}\label{31}
1=k_{n,0}^*\le k_{n,1}^*\le k_{n,2}^*\le k_{n,3}^*\le  k_{n,4}^*\le k_{n,5}^*=j,\qquad j\to\infty,\end{equation}
where positive integers $\{k_{n,i}^*\}$
are such that
\begin{equation}\label{4a}\{k_{n,i-1}^*=o(k_{n,i}^*), i\in\{1,2,...,5\}\}.
\end{equation}
Roughly speaking, the partition is required to split the conditional probability of the  maximum $M_{1,j}$ in $P\{T_1(u)=j\}=P\{M_{1,j}\le u,X_{j+1}>u|X_{1}>u\}$ into the product of independent probabilities of partial maxima $M_{1,k^*_{n,1}}$, $M_{k^*_{n,2},k^*_{n,3}}$ and $M_{k^*_{n,4},j}$. The independence follows from mixing conditions (\ref{29}), (\ref{29a}). The statement (\ref{2a}) is obtained from Theorem 2.1 and   Corollary 2.3 of \cite{obrien}.
\begin{theorem}\label{T2}
Let  $\{X_n\}_{n\ge 1}$ be a stationary
process with the extremal index
$\theta$.
Let $\{x_{\rho_n}\}$ be a sequence of quantiles of $X_1$ of the levels $\{1-\rho_n\}$,
that satisfies the conditions (\ref{1}) and (\ref{2}) if $u_n$ is replaced by $x_{\rho_n}$.
Let positive integers $\{k^*_{n,i}\}$, $i=\overline{0,5}$,
be   as in  (\ref{31}) and (\ref{4a}), respectively,
$\Delta_{n,i}=k_{n,i}^*-k_{n,i-1}^*$,
$q_{n,i}^*=o(\Delta_{n,i})$, $i\in\{1,2,...,5\}$,
 be such that
for each $\varepsilon>0$ there exist $n_{\varepsilon}$ and $j_0=j_0(n_{\varepsilon})$ such that for all $n>n_{\varepsilon}$ and $j>j_0(n_{\varepsilon})$
\begin{eqnarray}\label{29}
\alpha^*_n(x_{\rho_n})&=&
\max\{
\alpha_{k_{n,4}^*,q_{n,1}^*};
\alpha_{k_{n,3}^*,q_{n,2}^*}; \alpha_{\Delta_{n,3},q_{n,3}^*}; \alpha_{j+1-k_{n,2}^*,q_{n,4}^*};\nonumber
\\
&&
\alpha_{j+1-k_{n,1}^*,q_{n,5}^*}
\}<\varepsilon\end{eqnarray}
and
\begin{eqnarray}\label{29a}
\alpha_{j+1,k_{n,4}^*-k_{n,1}^*}/\rho_n &<& \varepsilon
\end{eqnarray}
hold,
 where $\alpha_{n,q}=\alpha_{n,q}(x_{\rho_n})$ is determined by (\ref{44}).
 Then for the same $n$ and $j$ it holds
\begin{eqnarray}\label{2a}
|P\{T_1(x_{\rho_n})=j\}/(\theta^2\rho_n(1-\rho_n)^{(j-1)\theta})-1|&<&\varepsilon.
\end{eqnarray}
\end{theorem}
The theorem implies that the probability $P\{T_1(x_{\rho_n})=j\}$ is close to the geometric form corrupted by extremal index $\theta$ for sufficiently large $n$ and $j$.
\\
The paper is organized as follows.
In Section \ref{Sec2} we derive the limit
distribution and 
expectation of the  first hitting time to exceed a sufficiently high threshold. The limit distribution of the second hitting time is obtained in Section \ref{Sec3}.  In Section \ref{Sec1} examples of first hitting time distributions are obtained for different processes including real data. Conclusions are given in Section \ref{sec-conclusion}. Proofs are presented in the Appendix.

\section{Distribution and expectation of the first hitting time}\label{Sec2}

For all $n$ and $j$ sufficiently large one can rewrite (\ref{2a}) in a geometric form as
\begin{eqnarray}
\label{10}
|\frac{c_nP\{T_1(x_{\rho_n})=j\}}{\eta_n(1-\eta_n)^{j-1}}-1|<\varepsilon,
\end{eqnarray}
where
$c_n=\eta_n/\left(\theta^2\left(1-(1-\eta_n)^{1/\theta}\right)\right)$,
$0< \eta_n < 1$,
using the replacement $(1-\rho_n)^{\theta}=1-\eta_n$.
We shall use (\ref{10}) to prove the next theorem.
\begin{theorem}\label{Theor1} Let all conditions of Theorem \ref{T2} be satisfied.
Then for the same $n$ and $j$ as in Theorem \ref{T2} we get
\begin{eqnarray}\label{19}|\frac{P\{T^*(x_{\rho_n})=j\}}{\psi_{j-1}(n)}-1|&<&\varepsilon,
\end{eqnarray}
where
\begin{eqnarray}\label{17} \psi_{j-1}(n)&=& \frac{\theta^2\rho_n^2(1-\rho_n)^{\theta (j-1)}}{1-(1-\rho_n)^{\theta}}.\end{eqnarray}
 \end{theorem}
From (\ref{1}) \begin{eqnarray}\label{025}\rho_n&\sim &\tau/n\qquad \mbox{and}\qquad (1-\rho_n)^{\theta}=1-\theta\rho_n+o(\rho_n)\end{eqnarray}  hold as  $n\to\infty$.
Expressions (\ref{19}) and (\ref{17})
imply that for any positive $\varepsilon$ there exists $n_{\varepsilon}$ such that for $n>n_{\varepsilon}$ and $j>j_0(n_{\varepsilon})$ the  probability of the first hitting time has a geometric distribution with probability $\theta\rho_n$, i.e.
\begin{eqnarray*}|\frac{P\{T^*(x_{\rho_n})=j\}}{\theta\rho_n(1-\theta\rho_n)^{j-1}}-1|&<&\varepsilon.
\end{eqnarray*}
Together with (\ref{2a}) it implies that for sufficiently large $n$ and $j$ it holds
\[P\{T^*(x_{\rho_n})=j\}\approx \theta P\{T_1(x_{\rho_n})=j\}.\]
\begin{lemma}\label{Lem1}Let the conditions of Theorem \ref{T2} be satisfied and for some $\beta>0$
\begin{equation}\label{024}
\sup_{n}E ((T^*(x_{\rho_n}))^{1+\beta})/\Lambda_{n}<\infty
\end{equation}
holds.
Then it follows
\begin{equation}\label{5a}
|ET^*_{j_0}(x_{\rho_n})/(\Lambda_n\rho_n)-1|<\varepsilon,
\end{equation}
where $j_0=o(n)$, $ET^*_{j_0}(x_{\rho_n})=\sum_{j=j_0+1}^{\infty}jP\{T^*(x_{\rho_n})=j\}$,
\begin{equation}\label{5b}\Lambda_n=\frac{\theta^2\rho_n}{(1-(1-\rho_n)^{\theta})^3}
(1-\rho_n)^{\theta j_0}\left(j_0(1-(1-\rho_n)^{\theta})+1\right).
\end{equation}
\end{lemma}
  The expression (\ref{5a}) specifies the rate of convergence in (\ref{5bc}).
\begin{remark}  The condition (\ref{024}) provides a uniform convergence of the series $\sum_{j=1}^{\infty}jP\{T^*(x_{\rho_n})=j\}/\Lambda_n$ by $n$. The condition is fulfilled particularly for $T^*(x_{\rho_n})$ corresponding to the ARMAX process, see Section \ref{Sec4.1}. From (\ref{8bb}) we have
$\sup_{n}E ((T^*(x_{\rho_n}))^{1+\beta})/\Lambda_{n}<\sup_{n}E ((T^*(x_{\rho_n}))^{2})/\Lambda_{n}\sim \sup_{n}(2-\theta\rho_n)(1-\rho_n)^{1-\theta(j_0+1)}/(\theta(1+j_0\theta\rho_n))<\infty$ for $\rho_n\sim \tau/n$.
\end{remark}
Let us turn to (\ref{10a}). If $\{X_i\}$ and $\{Y_i\}$ are mutually independent then
\begin{eqnarray*}P\{T_T^*(u)=j+1\}
&=& P\{M_j\le u, X_{j+1}>u\}P\{S_{j+1}\le T\}
\\
&=& P\{T^*(u)=j+1\}P\{\sum_{i=1}^jY_i\le T\}
\end{eqnarray*}
follows. From (\ref{19}) and (\ref{17}) we get
\begin{eqnarray*}|P\{T_T^*(x_{\rho_n})=j+1\}/(\psi_j(n)(1-P\{\sum_{i=1}^jY_i> T\}))-1|<\varepsilon
\end{eqnarray*}
for any $\varepsilon>0$ and $n>n_{\varepsilon}$ and $j>j_0(n_{\varepsilon})$.
Assuming $Y_i$'s are iid regularly varying random variables with tail index $\alpha\ge 0$  we have
\begin{eqnarray*}P\{\sum_{i=1}^jY_i> T\}&\sim & jP\{Y_i> T\}\sim jT^{-\alpha}
\end{eqnarray*} for $j\ge 1$ as $T\to\infty$, see
Lemma 3.1, \cite{jessen}.
The condition $1-jT^{-\alpha}>0$ is provided by $T>j_0^{1/\alpha}$ since $j>j_0$.
Then the lemma follows.
\begin{lemma}\label{Lem2} Let the conditions of Theorem \ref{T2}
be satisfied. Let $\{X_i\}$ and $\{Y_i\}$ in (\ref{10a}) be mutually independent and  $\{Y_i\}_{i\ge 1}$ be iid regularly varying random variables with tail index $\alpha\ge 0$ and $T>j_0^{1/\alpha}$ holds. Then for the same $n$ and $j$ as in Theorem \ref{T2} we get
 \begin{eqnarray*}|P\{T_T^*(x_{\rho_n})=j+1\}/(\psi_j(n)(1-jT^{-\alpha})-1|<\varepsilon.
 \end{eqnarray*}
\end{lemma}

\section{Distribution of the second hitting time
}\label{Sec3}
Let us denote the second hitting time of $u_n$ as $T^{**}(u_n)$.
The probability to hit $u_n$ twice is determined by
      \begin{eqnarray}\label{15} && P\{T^{*}(x_{\rho_n})=j, T^{**}(x_{\rho_n})=j+m\}
      \\
      &=& P\{M_{j-1}\le u_n, X_j>u_n, M_{j,j+m-1}\le u_n, X_{j+m}>u_n\},\qquad m=1, 2,...\nonumber
      \end{eqnarray}
      \begin{lemma}\label{Lem3} Let the conditions of Theorem \ref{T2} be satisfied. Then for the same $n$ and $j$ as in Theorem \ref{T2}
       we have
      \begin{eqnarray*}
      |\frac{P\{T^{*}(u_n)=j, T^{**}(u_n)=j+m\}}{P\{\chi=j\}P\{\chi=m\}}-1|<\varepsilon,
      \end{eqnarray*}
       where $\chi$ is a geometrically distributed random variable  with probability $\rho_n\theta$.
           \end{lemma}
      Similarly the statement can be extended to the probability of the $k$th hitting time, i.e. the minimal time to find $k$ large nodes of the network.  Random walks used in social networks as sampling may return to the same nodes with some positive probability. This may reduce the number of distinct nodes in the sample and particularly ones which degrees exceed the threshold.  The degrees of repeated nodes may not exceed the threshold and hence, do not impact on the probability to reach $k$ different large nodes. Moreover,  the degrees of repeated nodes may change over time. These problems are out of scope of this paper.
\section{Examples}\label{Sec1}
\subsection{ARMAX process}\label{Sec4.1}
Let us obtain the distribution of the first hitting time of the ARMAX process. The latter process is determined as \[X_t=\max\{\alpha X_{t-1},(1-\alpha)Z_t\}, \qquad t\in \textmd{Z},\]
where $0\le \alpha < 1$, and  $\{Z_t\}$ are iid standard Fr\'{e}chet distributed r.v.s with the distribution function $F(x)=\exp\left(-1/x\right)$, $x>0$. The r.v. $X_t$ has the same distribution assuming $X_0=Z_0$. The extremal index is equal to $\theta=1-\alpha$, \cite{ancona00}.
\\
Using that \begin{eqnarray}\label{022}P\{X_i\le x_{\rho}\}=1-\rho=q=e^{-1/x_{\rho}}\end{eqnarray} and
\begin{eqnarray}\label{023}P\{\alpha Z_i\le x_{\rho}\}&=& e^{-\alpha/x_{\rho}}=(1-\rho)^{\alpha}=q^{\alpha}\end{eqnarray}
 we derive in Section \ref{ProofARMAX} the following.
 \begin{proposition}\label{P1} For the ARMAX and MM processes we have
\begin{eqnarray}\label{8c}P\{T^*(x_{\rho})=j\}&=&(1-(1-\rho)^{\theta})(1-\rho)^{\theta(j-2)+1},
\end{eqnarray}
\begin{eqnarray}\label{8a}ET^*(x_{\rho})&=&(1-\rho)^{1-\theta}/(1-(1-\rho)^{\theta})
\end{eqnarray}
and
\begin{eqnarray}\label{8bb}E(T^*(x_{\rho}))^2&=&(1-\rho)^{1-\theta}(1+(1-\rho)^{\theta})/(1-(1-\rho)^{\theta})^2.
\end{eqnarray}
\end{proposition}

\subsection{MM process}
\begin{figure*}[th!]
\caption{The distribution (\ref{8c}) of the first hitting time of the ARMAX and MM processes
and the model (\ref{17}) with $\theta=0.1$
for $j=5$ (top) and  $j=20$ (bottom) against $\rho$, where $\rho$ close to zero corresponds to a high quantile $x_{\rho}$ as the threshold $u$.}
\label{fig:1}
\centerline{\epsffile{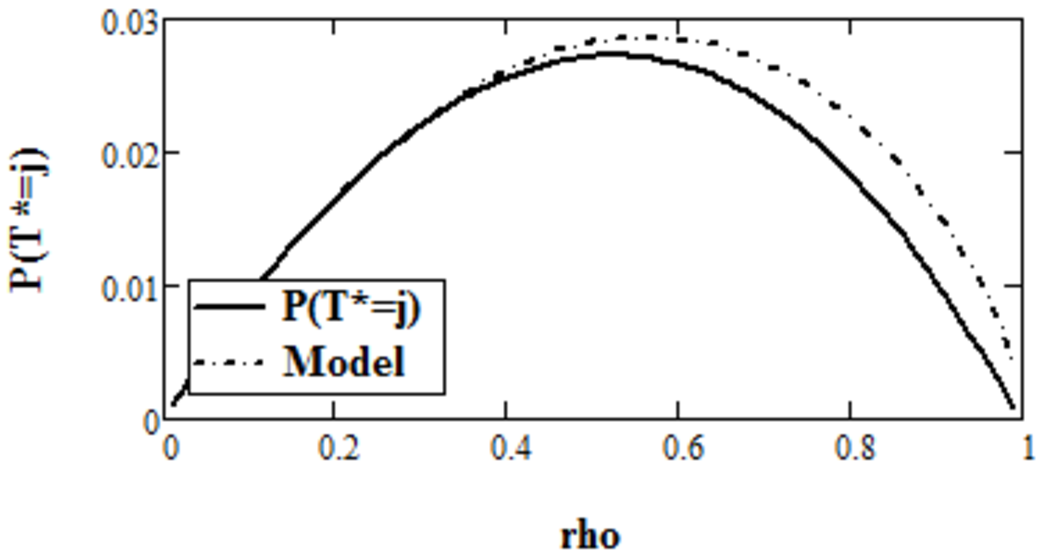}}

\centerline{\epsffile{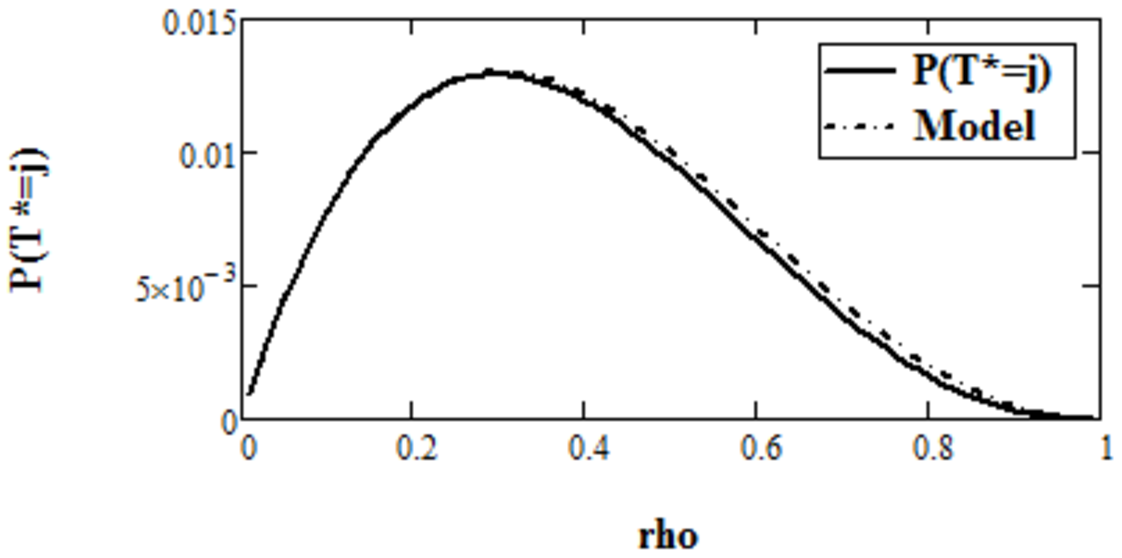}}
\end{figure*}

\begin{figure*}[t!]
\caption{The mean first hitting time (\ref{8a}) of the ARMAX and MM processes
and the model $\Lambda_n\rho_n$ based on (\ref{5a}) and (\ref{5b}) with $\theta=0.1$
for $j_0=0$ and  $j_0=5$ against $\rho$.}
\label{fig:2}
\centerline{\epsffile{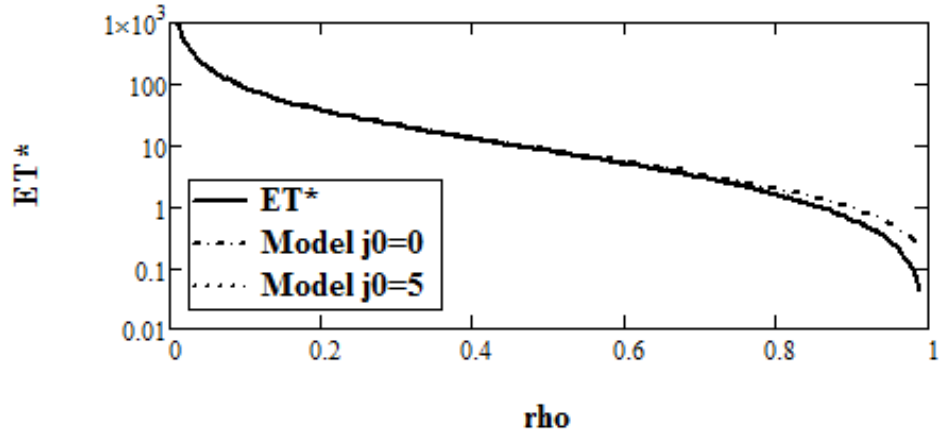}}
\end{figure*}

We obtain the distribution of the first hitting time of the MM process. This process is determined by the formula
\[X_t=\max_{i=0,...,m}\{\alpha_iZ_{t-i}\}, \qquad t\in\textmd{Z},\]
where  $\{\alpha_i\}$ are nonnegative constants such that $\sum_{i=0}^m \alpha_i=1$ and   $\{Z_t\}$ are iid standard Fr\'{e}chet distributed r.v.s.
 The distribution of $X_t$ is also standard Fr\'{e}chet. The extremal index of the process is determined by $\theta=\max_i\{\alpha_i\}$, \cite{ancona00}.
  Assuming $\alpha_0\ge \alpha_1\ge ...\ge \alpha_m$ we derive in Section \ref{ProofMM} that the distribution of the first hitting time is the same as for the ARMAX process.
\\
In Figure \ref{fig:1} the comparison of the exact distribution of $T^*(x_{\rho})$ for the ARMAX and MM processes and the model obtained in Theorem \ref{Theor1} is shown. The model (\ref{17}) is valid for sufficiently large $n$ and $j$. This corresponds to $\rho_n$ close to zero and high quantiles $x_{\rho_n}$ using as thresholds $u_n$. Thus the model approximates the distribution (\ref{8c}) better for small $\rho$ and large $j$.
\\
The comparison of the mean first hitting time (\ref{8a})  and the theoretical model obtained in Lemma \ref{Lem1} is shown in Figure \ref{fig:2}. The difference is observed only for $j_0=0$ and when $\rho$ is close to $1$. It should be noted that we consider $\theta=0.1$ corresponding to a large local  dependence in the extremes of the process $\{X_n\}$. $\theta=1$ corresponds to independent observations.
\subsection{AR(1) process}
We consider the AR(1) process with uniform noise, \cite{chernick}. For a fixed integer $r\geq 2$ let $\epsilon_n$, $n\geq 1$ be iid r.v.s with
$P\{\epsilon_1=k/r\}=1/r$, $k=0,1, \ldots , r-1$.
The process is defined by \begin{eqnarray*}\label{11a}X_j&=&(1/r) X_{j-1}+\epsilon_j, \quad j \geq 1\qquad \mbox{and} \quad X_0\sim  U(0,1)\end{eqnarray*}
with $X_0$ independent of the $\epsilon_j$. Since $X_0\sim  U(0,1)$ then $X_1\sim  U(0,1)$ holds. The extremal index of AR(1) is $\theta=1-1/r$.
\begin{proposition}\label{P2} For the AR(1) process we have
\begin{eqnarray}\label{8}P\{T^*(u_n)=j\}&=&\left\{
                     \begin{array}{lll}
                       1-\theta, & j=1 \\
                       (1-\theta)^j\left(u_n-j\theta(1-u_n)\right),
                       & 2\leq j \leq j_0 \\
                       (1-\theta)^{j_0+2}\left(u_n-j\theta(1-u_n)\right),
                       & j_0<j\le m-1,
                     \end{array}
                   \right.
\end{eqnarray}
where $j_0=[\ln n/(2\ln r)]$ and $m$ satisfies the inequality
\begin{eqnarray}\label{7aa}&&-\frac{\ln(1-u_n)}{\ln(r)}-1<m-1\le -\frac{\ln(1-u_n)}{\ln(r)}.\end{eqnarray}
\end{proposition}
Selecting $u_n=1-x/n$, $x>0$ we get $u_n-j\theta(1-u_n)=1-(x/n)(1+j\theta)$ that is positive for sufficiently large $n$.
\\
The proof is given in Section \ref{ProofAR}.
\begin{remark} The mixing conditions (\ref{29}) and (\ref{29a}) of Theorem \ref{T2} are fulfilled for the ARMAX and the AR(1) processes if $j>j_0(n)$ holds, where $j_0(n)\to\infty$ as $n\to\infty$, and for the MM process if $j>m$ and $\alpha_0\ge \alpha_1\ge...\ge \alpha_m$ hold, \cite{markovich16b}.
\end{remark}
\subsection{Real data}
\begin{figure*}[t!]
\caption{The model $\Lambda_n\rho_n$ of the mean first hitting time  calculated by (\ref{5a}) and (\ref{5b}) of the Enron and DBPL data sets
 with $\theta=0.22$ and $\theta=0.15$, respectively,
for $\rho=0.05$  against $j_0$.}
\label{fig:3}
\centerline{\epsffile{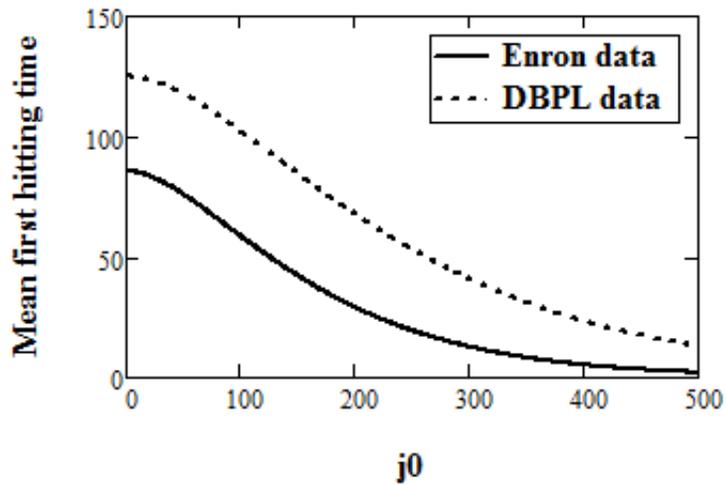}}
\end{figure*}

\begin{figure*}[t!]
\caption{The sample autocorrelation function (\ref{026}) of the Enron and DBPL data sets.}
\label{fig:4}
\centerline{\epsffile{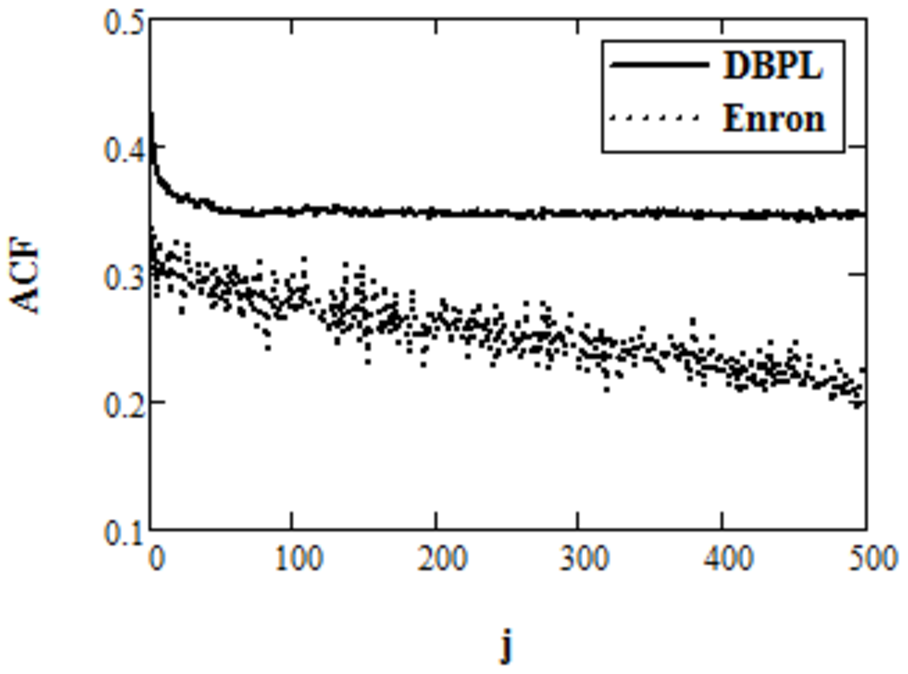}}
\end{figure*}

We consider two real data sets of the Enron email and DBPL networks presented in \cite{snapnets} and investigated in \cite{markovich15}. The sets contain node degrees. In \cite{markovich15} it was found that the both data sets are heavy-tailed distributed and their extremal index $\theta$ was calculated by intervals estimator proposed by \cite{ferro}. Typically, it may be assumed that the node degrees are regularly varying distributed.  In Figure \ref{fig:3} the model from Lemma \ref{Lem1} of the mean first hitting time against $j_0$ is shown. The $j_0$ indicates the truncated expectation $ET^*_{j_0}(x_{\rho})$. The theoretical model is valid for any sampling technique (a random walk, Markov chain) that satisfies the mixing condition  (\ref{29})-(\ref{29a})  for $j>j_0$, where  $j_0=j_0(n)$ is sufficiently large. The latter condition is equivalent to the $j-$dependence.   The $j-$dependence may be checked in practice by an autocorrelation function (ACF) (see, e.g., \cite{markovichkrieger}). Since both data sets have infinite variance according to \cite{markovich15} it is better to use the special sample ACF for heavy-tailed data recommended in \cite{davisresnick}, i.e.
\begin{equation}\label{026}\widetilde{\rho}(j)=\sum_{t+1}^{n-j}X_tX_{t+j}/\sum_{t=1}^nX_t^2\end{equation} at lag $j$. This ACF is not centralized by the sample average $\overline{X}$ in contrast to the classical sample ACF. Moreover, this estimate may behave in a very unpredictable way
if one uses the class of non-linear processes in the sense that
$\widetilde{\rho}(j)$ may converge in distribution to a
non-degenerate random variable depending on $j$.
For linear processes it converges in distribution to a constant
depending on $j$, \cite{davisresnick}. From Figure \ref{fig:4} one may conclude that the DBPL data are short-range dependent since its ACF decreases after $j\approx 50$ as far as the Enron data are not. This may indirectly indicate that the DBPL and Enron data determine linear and  non-linear processes, respectively.
\\
Everything what we need for our model are the extremal index $\theta$ and the quantile threshold $1-\rho$. We take $\rho=0.05$ that corresponds to $95\%$ quantile $x_{\rho}$ of an underlying data set taken as the threshold. We may conclude from Figure \ref{fig:3} that the mean minimal time required to reach a  node with  degree larger than $u=x_{\rho}$ is longer for the DBPL data than for the Enron data.
\section{Conclusions}\label{sec-conclusion} We have obtained the limit distribution and expectation of the first hitting time for processes which satisfy the mixing conditions of Theorem \ref{T2}. The latter are fulfilled particularly for Markov chains represented by the ARMAX, the MM and the AR(1) processes. Exact distributions of the first hitting time for the latter processes are obtained. Markov chains used in social networks as sampling techniques can be compared with regard to the quantiles and expectation of the first hitting time. The presented research can be particularly useful for such comparison of sampling strategies. The results are extended to the second hitting time.

\section*{Acknowledgements}
I appreciate the financial support by the  DAAD scholarship 57210259 (Germany) 2016.

\section{Appendix}\label{Proof}

\subsection{Proof of Theorem \ref{Theor1}}
It follows from (\ref{0}) that
\begin{eqnarray}\label{3}P\{T^*(u_n)=j+1\} &=& P\{M_{j}\le u_n, X_{j+1}>u_n\}
\\
&=& P\{M_{j}\le u_n\}-P\{M_{j+1}\le u_n\}.\nonumber
\end{eqnarray}
Following Ferro and Segers (2003)  we get alternatively for $n\ge 1$
\begin{eqnarray*}P\{T_1(u_n)>j\} &=& P\{M_{1,j+1}\le u_n| X_{1}>u_n\}
\\
&=&\left(P\{M_{1,j+1}\le u_n\}-P\{M_{j+1}\le u_n\}\right)/P\{X_{1}>u_n\}
\\
&=&\left(P\{M_{j}\le u_n\}-P\{M_{j+1}\le u_n\}\right)/P\{X_{1}>u_n\}
\\
&=& P\{T^*(u_n)=j+1\}/P\{X_{1}>u_n\}.
\end{eqnarray*}
Thus, we get
\begin{eqnarray*}\label{4}
P\{T^*(x_{\rho_n})=j+1\}&=& P\{X_{1}>x_{\rho_n}\}\cdot P\{T_1(x_{\rho_n})>j\}\nonumber
\\
&=& \rho_n\sum_{i=j+1}^{\infty}P\{T_1(x_{\rho_n})=i\}.
\end{eqnarray*}
From (\ref{10})
 we obtain
\begin{eqnarray}\label{5}
P\{T^*(x_{\rho_n})=j+1\}&=& \frac{\rho_n}{c_n}\sum_{i=j+1}^{\infty}c_nP\{T_1(u_n)=i\}
\\
&<& (1+\varepsilon)\frac{\rho_n}{c_n}\sum_{i=j+1}^{\infty}\eta_n(1-\eta_n)^{i-1}=(1+\varepsilon)\psi_{j-1}(n),\nonumber
\\
P\{T^*(x_{\rho_n})=j+1\}&>&  (1-\varepsilon)\psi_{j-1}(n),\nonumber
\end{eqnarray}
where $\psi_{j-1}(n)$ is determined by (\ref{17}).
Since $\psi_{j-1}(n)\to 0$ as $n\to\infty$ holds, the series in  (\ref{5}) converges uniformly by all $n>n_{\varepsilon}$ by Weierstrass' theorem.

\subsection{Proof of Lemma \ref{Lem1}}
 Let us consider
 the expectation of the first hitting time
\begin{eqnarray*}
ET^*(x_{\rho_n})&=& \sum_{j=1}^{\infty}jP\{T^*(x_{\rho_n})=j\}.
\end{eqnarray*}
From (\ref{17}) we get
\begin{eqnarray*}
&&\sum_{j=j_0+1}^{\infty}j\psi_{j-1}(n)=\sum_{j=j_0+1}^{\infty}j\frac{\theta^2\rho_n^2(1-\rho_n)^{\theta(j-1)}}{1-(1-\rho_n)^{\theta}}=\Lambda_n\rho_n,
\end{eqnarray*}
where $\Lambda_n$ is determined by (\ref{5b}).
Due to (\ref{025}) we have
\[\Lambda_n\sim (1-\theta\rho_n)^{j_0}\left(j_0\theta\rho_n+1\right)/(\theta\rho_n^2) \sim \exp(-\tau\theta j_0/n)/(\theta\rho_n^2)
\to\infty\]
as $n\to\infty$.
Let us denote
\begin{eqnarray*}\label{33}
a_j(n)=jP\{T^*(x_{\rho_n})=j\}\end{eqnarray*}
and
\begin{eqnarray*}\label{35}
S_k(n)=\sum_{j=1}^{k-1}a_j(n)/\Lambda_{n},
\qquad r_k(n)=\sum_{j=k}^{\infty}a_j(n)/\Lambda_{n}.\end{eqnarray*}
We have to prove that $S(n)=\sum_{j=1}^{\infty}a_j(n)/\Lambda_{n}$
converges uniformly by $n$. 
For this purpose, we shall prove that \[\lim_{k\to\infty}\sup_{n} r_k(n)=0.\]
The latter follows  from
\begin{eqnarray}\label{11}
\sup_{n}r_k(n)&=& \sup_{n}\sum_{j=k}^{\infty}\frac{jk^{\beta}
P\{T^*(x_{\rho_n})=j\}}{k^{\beta}\Lambda_{n}}
\le
 \frac{1}{k^{\beta}}
 \sup_{n} \frac{E(T^*(x_{\rho_n}))^{1+\beta}}{\Lambda_{n}}\nonumber
\end{eqnarray}
and the assumption (\ref{024}).
It remains to prove that
\begin{equation}\label{34}
\lim_{n\to\infty}S_{j_0}(n)/\rho_n=1,
\end{equation}
were $S_{j_0}(n)=\sum_{j=j_0+1}^{\infty}a_j(n)/\Lambda_{n}$.
\\
Using the replacement $(1-\rho_n)^{\theta}=1-\eta_n$ from  (\ref{19}), (\ref{17}) and (\ref{5b}) we get for any $\varepsilon>0$ that it holds
\begin{eqnarray*}
S_{j_0}(n)&<&  \frac{(1+\varepsilon)(1-(1-\rho_n)^{\theta})^2}{(1-\rho_n)^{\theta j_0}(j_0(1-(1-\rho_n)^{\theta})+1)}\sum_{j=j_0+1}^{\infty}
j\rho_n(1-\rho_n)^{(j-1)\theta}
\\
&=&
\frac{(1+\varepsilon)\eta_n\rho_n}{(1-\eta_n)^{j_0}(j_0\eta_n+1)}
\sum_{j=j_0+1}^{\infty}
j\eta_n(1-\eta_n)^{j-1}.
\end{eqnarray*}
Similarly, one can get
\begin{eqnarray*}
S_{j_0}(n)&>&\frac{(1-\varepsilon)\eta_n\rho_n}{(1-\eta_n)^{j_0}(j_0\eta_n+1)}
\sum_{j=j_0+1}^{\infty}
j\eta_n(1-\eta_n)^{j-1}.\end{eqnarray*}
Since 
\begin{eqnarray*}
&&\sum_{j=j_0+1}^{\infty}j\eta_n(1-\eta_n)^{j-1}=\frac{(1-\eta_n)^{j_0}}{\eta_n}(j_0\eta_n +1)\end{eqnarray*}
 and $\varepsilon$ is arbitrary then (\ref{34}) and thus,
(\ref{5a}) follows.
\subsection{Proof of Lemma \ref{Lem3}}
 By (\ref{18})  and the stationarity of $\{X_n\}$  we obtain
\begin{eqnarray}\label{13}&& P\{T_1(u_n)=n\} = P\{M_{1,n}\le u_n, X_{n+1}>u_n|X_{1}>u_n\}\nonumber
\\
&=& \left(P\{M_{1,n}\le u_n, X_{n+1}>u_n\}-P\{M_{n}\le u_n, X_{n+1}>u_n\}\right)/P\{X_{1}>u_n\}\nonumber
\\
&=& \left(P\{M_{n-1}\le u_n, X_{n+1}>u_n\}-P\{M_{n}\le u_n, X_{n+1}>u_n\}\right)/P\{X_{1}>u_n\}\nonumber
\\
&=&\left(P\{T^*(u_n)=n\}-P\{T^*(u_n)=n+1\}\right)/P\{X_{1}>u_n\}.
\end{eqnarray}
 From (\ref{15}) we get due to stationarity
      \begin{eqnarray*} && P\{T^{*}(u_n)=j, T^{**}(u_n)=j+m\}
      \\
      &=& P\{M_{j-1}\le u_n, M_{j,j+m-1}\le u_n, X_{j+m}>u_n\}-P\{M_{j+m-1}\le u_n, X_{j+m}>u_n\}
      \\
      &=& P\{M_{j+m-2}\le u_n, X_{j+m-1}>u_n\}-P\{T^{*}(u_n)=j+m\}
      \\
      &=& P\{T^{*}(u_n)=j+m-1\}-P\{T^{*}(u_n)=j+m\}
      \\
      &=& P\{T_1(u_n)=j+m-1\}P\{X_1>u_n\}.
      \end{eqnarray*}
      The last two lines are obtained from (\ref{0}) and (\ref{13}).
      Then using (\ref{10}) and denoting $\varphi_{j+m-2}(n)=\rho_n\eta_n(1-\eta_n)^{j+m-2}/c_n$ one can rewrite
      \begin{eqnarray*} && |P\{T^{*}(x_{\rho_n})=j, T^{**}(x_{\rho_n})=j+m\}/\varphi_{j+m-2}(n)-1|
      <\varepsilon.
      \end{eqnarray*}
      Since it holds
      \begin{eqnarray*}\varphi_{j+m-2}(n)&\sim &\theta^2\rho_n^2(1-\rho_n)^{(j+m-2)\theta}\sim\theta\rho_n(1-\theta\rho_n)^{j-1}\theta\rho_n(1-\theta\rho_n)^{m-1},
      \end{eqnarray*}
      the statement of the lemma follows.

\subsection{Proof of Proposition \ref{P1} for an ARMAX process}\label{ProofARMAX}
From the definition of the ARMAX process we obtain the distribution of the first hitting time. It holds
\begin{eqnarray*}&&P\{T^*(u)=j\}=P\{M_{j-1}\le u, X_j>u\}
\\
&=& P\{X_1\le u,..,X_{j-1}\le u, X_j>u\}
\\
&=&P\{X_1\le u, (1-\alpha)Z_2\le u,...,(1-\alpha)Z_{j-1}\le u, \max\{\alpha X_{j-1}, (1-\alpha)Z_j\}>u\},
\end{eqnarray*}
since $X_{i+1}\le u$, $i=1,...,j-2$ leads to $\alpha X_{i}\le u$ and $(1-\alpha)Z_{i+1}\le u$ and together with $X_i\le u$ it implies both $X_i\le u$ and $(1-\alpha)Z_{i+1}\le u$ due to $0<\alpha<1$. For an independent sequence $\{X_t\}$   $\alpha=0$ holds.
\\
Let us consider $\max\{\alpha X_{j-1}, (1-\alpha)Z_j\}>u$. Supposing $\alpha X_{j-1}>u$ contradicts $X_{j-1}\le u$. Hence, it follows $(1-\alpha)Z_j>u$ and it holds
\begin{eqnarray*}P\{T^*(u)=j\}&=&P\{X_1\le u, (1-\alpha)Z_2\le u,...,(1-\alpha)Z_{j-1}\le u, (1-\alpha)Z_j>u\}.
\end{eqnarray*}
Taking the $(1-\rho)$-level quantile $x_{\rho}$ as $u$ and using (\ref{022}) and (\ref{023}) we get
\begin{eqnarray*}P\{T^*(x_{\rho})=j\}&=& (1-\rho)(1-\rho)^{(1-\alpha)(j-2)}(1-(1-\rho)^{1-\alpha})
\end{eqnarray*}
and (\ref{8c}) follows.
\\
We shall obtain $ET^*(x_{\rho})$ for the ARMAX process. Denoting $1-\eta=(1-\rho)^{\theta}$ we get
\begin{eqnarray*}ET^*(x_{\rho})&=&\sum_{j=1}^{\infty}jP\{T^*(x_{\rho})=j\}=(1-\rho)(1-(1-\rho)^{\theta})\sum_{j=1}^{\infty}j(1-\rho)^{\theta(j-2)}
\\
&=& (1-\eta)^{1/\theta-1}\sum_{j=1}^{\infty}j\eta(1-\eta)^{j-1}=(1-\eta)^{1/\theta}/(\eta(1-\eta)).
\end{eqnarray*}
Similarly, one can get (\ref{8bb}).
\subsection{Proof of  Proposition \ref{P1} for a MM process}\label{ProofMM}
From the definition of the MM process we get
\begin{eqnarray}\label{9}&&P\{T^*(u)=j\}
\\
&=& P\{\max_{i=0,...,m}\{\alpha_iZ_{1-i}\}\le u,..,\max_{i=0,...,m}\{\alpha_iZ_{j-1-i}\}\le u, \max_{i=0,...,m}\{\alpha_iZ_{j-i}\}>u\}\nonumber
\end{eqnarray}
Assuming $\alpha_0\ge \alpha_1\ge ...\ge \alpha_m$ we obtain that the right-hand side of (\ref{9}) is equal to
\begin{eqnarray*}&&P\{\alpha_mZ_{1-m}\le u,.., \alpha_0Z_{1}\le u, \alpha_0Z_{2}\le u,..., \alpha_0Z_{j-1}\le u, \max_{i=0,...,m}\{\alpha_iZ_{j-i}\}>u\}.
\end{eqnarray*}
Let us consider the event $\{\max_{i=0,...,m}\{\alpha_iZ_{j-i}\}>u\}$. This is equivalent to $\{\alpha_0Z_{j}> u, \alpha_1Z_{j-1}\le u,..., \alpha_mZ_{j-m}\le u\}$.
\\
Really, suppose $\alpha_mZ_{j-m}> u$ holds. But this is in contradiction with $\alpha_{m-1}Z_{j-m}\le u$ in (\ref{9}). Furthermore, $\alpha_1Z_{j-1}> u$ contradicts $\alpha_0Z_{j-1}\le u$ etc.
\\
Summarizing we obtain
\begin{eqnarray*}&&P\{T^*(u)=j\}
= P\{\alpha_mZ_{1-m}\le u,.., \alpha_0Z_{1}\le u, \alpha_0Z_{2}\le u,..., \alpha_0Z_{j-1}\le u,
\\
&&\qquad \qquad \qquad \qquad\qquad\qquad\qquad\alpha_mZ_{j-m}\le u,..., \alpha_1Z_{j-1}\le u, \alpha_0Z_{j}> u\}
\\
&=& P\{\alpha_mZ_{1-m}\le u,.., \alpha_0Z_{1}\le u, \alpha_0Z_{2}\le u,..., \alpha_0Z_{j-1}\le u, \alpha_0Z_{j}> u\}.
\end{eqnarray*}
Hence, from (\ref{022}) and (\ref{023}) it follows
\begin{eqnarray*}P\{T^*x_{\rho})=j\}&=& q^{\alpha_m+...+\alpha_0+ \alpha_0(j-2)}(1-q^{\alpha_0})=(1-(1-\rho)^{\alpha_0})(1-\rho)^{1+\alpha_0(j-2)}.
\end{eqnarray*}
Thus, (\ref{8c}) follows.
\subsection{Proof of  Proposition \ref{P2}}\label{ProofAR}
To prove (\ref{8}) we use (\ref{3}) and results obtained in \cite{chernick1981} and \cite{chernick}.
For $j=1$ we have
\begin{eqnarray*}P\{T^*(u_n)=j\} &=& P\{M_0\le u_n, X_1>u_n\}=P\{\epsilon_1=(r-1)/r\}=1-\theta,
\end{eqnarray*}
since $X_1>u_n$ implies $\epsilon_1=(r-1)/r$ for sufficiently large $n$ and $r<n/x$. Really, suppose $\epsilon_1\le (r-2)/r$ holds. Then we get $X_1=1/rX_0+\epsilon_1\le 1-1/r$. This contradicts to $X_{1}>u_n=1-x/n$ for $r<n/x$.
\\
From Lemma 2.5 in \cite{chernick} it follows that the event  $X_{j+1}>u_n=1-x/n$, $x>0$ with $1\le j\le j_0$ and $j_0=[\ln n/(2\ln r)]$ leads to
 \begin{eqnarray}\label{12a}\epsilon_2=\epsilon_3=...=\epsilon_{j+1}=(r-1)/r\end{eqnarray}
for all $n$ sufficiently large. From another side, from Lemma 2.6 in Chernick et al. (1991) it follows that for all $n$ sufficiently large, the event
$X_{j+1}>u_n$ for $j>j_0$ leads to
\begin{eqnarray}\label{3a}\epsilon_t=(r-1)/r, \quad t=j-j_0,...,~j+1.\end{eqnarray}
If $(r-1)x<n$ holds, we get by formula (4.3) in \cite{chernick1981}
\begin{eqnarray}\label{8b}P\{M_{j}\le u_n\}&=& 1-\frac{(j+1)r-j}{r}(1-u_n)\end{eqnarray}
and if $j\le m-1$ holds,
where $m$ is the integer for which $1-r^m(1-u_n)<0$ and $1-r^{m-1}(1-u_n)\ge 0$ (i.e. (\ref{7aa})) hold.
\\
 Thus, (\ref{3}) can be rewritten as
\begin{eqnarray*}P\{T^*(u_n)=j+1\} &=& P\{M_{j}\le u_n\}P\{X_{j+1}>u_n\}.
\end{eqnarray*}
From (\ref{12a}) and (\ref{3a}) we get
\begin{eqnarray*}P\{X_{j+1}>u_n\} &=& \left\{
                     \begin{array}{ll}
                                         P\{\epsilon_t=(r-1)/r,~ t=2,...,~j+1\},
                       & \mbox{if}~ 2\leq j \leq j_0, \\
                       P\{\epsilon_t=(r-1)/r, ~ t=j-j_0,...,~j+1\},
                       & \mbox{if}~ j_0<j\le m-1
                     \end{array}
                   \right.
\\
&=&\left\{
                     \begin{array}{ll}
                                              \prod_{t=2}^{j+1}P\{\epsilon_t=(r-1)/r\},
                       & \mbox{if}\qquad  2\leq j \leq j_0, \\
                       \prod_{t=j-j_0}^{j+1}P\{\epsilon_t=(r-1)/r\},
                       & \mbox{if}\qquad j_0<j\le m-1
                     \end{array}
                   \right.
\end{eqnarray*}
Then the statement follows from (\ref{8b}).


\begin{thebibliography}{99}

\bibitem[\protect\citeauthoryear{Ancona-Navarrete and Tawn}{2000}]{ancona00}
Ancona-Navarrete, M.A.,  Tawn, J. A. (2000) 'A comparison of Methods for Estimating the Extremal Index',
   {\it Extremes}, Vol. 3, No. 1, pp. 5--38.

\bibitem[\protect\citeauthoryear{Asmussen}{2000}]{asmussen1998}
Asmussen, S. (1998) 'Subexponential asymtotics for stochastic processes: extremal behavior, stationary distributions and first passage probabilities',
   {\it Ann. Appl. Probab.}, Vol. 8,  pp. 354--374.

\bibitem[\protect\citeauthoryear{Avrachenkov et al.}{2012}]{avrachenkov12}
Avrachenkov, K., Litvak, N., Sokol, M. and Towsley, D. (2012)
'Algorithms and Models for the Web Graph',   {\it Lecture Notes in Computer Science},  7323,
Quick Detection of Nodes with Large Degrees, Springer Berlin Heidelberg, pp. 54--65.

\bibitem[\protect\citeauthoryear{Avrachenkov et al.}{2015}]{avrachenkov15}
Avrachenkov K.,  Markovich N. M.,  Sreedharan J.K. (2015) 'Distribution and Dependence of Extremes in Network Sampling Processes',  {\it Computational Social Networks}, Vol. 2, No. 12, pp. 1--21,  doi:10.1186/s40649-015-0018-3

\bibitem[\protect\citeauthoryear{Beirlant et al.}{2004}]{beirlant}
Beirlant, J., Goegebeur, Y., Teugels, J. and Segers, J. (2004)
'Statistics of Extremes: Theory and Applications', {\it Wiley,
Chichester, West Sussex}

\bibitem[\protect\citeauthoryear{Bortot and Tawn}{1998}]{bortot}
Bortot, P. and Tawn, J.A. (1998) 'Models for the extremes of Markov chains', {\it Biometrika},   Vol. 85, pp. 851--867.


\bibitem[\protect\citeauthoryear{Chernick}{1981}]{chernick1981}
Chernick, M.R. (1981) 'A limit theorem for the maximum of autoregressive processes with uniform marginal distributions',
    {\it Ann. Probab.}, Vol. 9, pp. 145--149.

\bibitem[\protect\citeauthoryear{Chernick et al.}{1991}]{chernick}
Chernick, M.R.,   Hsing, T.  and  McCormick,  W.P. (1991)
'Calculating the extremal index for a class of stationary sequences',
   {\it  Adv. Appl. Prob.},  Vol. 23, pp. 835--850.

\bibitem[\protect\citeauthoryear{Choudhury and Medhi}{2011}]{Choudhury}
Choudhury, A. and Medhi, P. (2011) 'A simple analysis of customer impatience in multi-server queues', {\it Int. J. of Applied Management Science}, Vol. 3, No. 3,  pp. 294-–315.

\bibitem[\protect\citeauthoryear{Lee et al.}{2012}]{lee}
 Chul-Ho Lee, Xin Xu and Do Young Eun. (2012)	'Beyond Random Walk and Metropolis-Hastings Samplers: Why You Should Not Backtrack for Unbiased Graph Sampling', {\it CoRR}.

\bibitem[\protect\citeauthoryear{D\'{a}n, G. and  Fodor}{2009}]{dan}
  D\'{a}n, G. and  Fodor, V. (2009) 'Delay asymptotics and scalability for peer-to-peer live streaming', {\it IEEE Trans. Parallel Distrib.}, Vol. 20, No. 10,  pp. 1499-–1511.

\bibitem[\protect\citeauthoryear{Davis and Resnick}{1985}]{davisresnick}
Davis, R. and Resnick, S. (1985) 'Limit theory for moving
averages of random variables with regularly varying tail
probabilities', {\it Ann. Probability}, Vol. 13, pp. 179--195.

\bibitem[\protect\citeauthoryear{Doukhan et al.}{2015}]{doukhan}
Doukhan, P., Jakubowski, A. and Lang, G. (2015) 'Phantom distribution functions for some stationary
sequences', {\it Extremes}, Vol. 18,  pp. 697–-725.

\bibitem[\protect\citeauthoryear{Ferro and  Segers}{2003}]{ferro}
  Ferro, C.A.T. and  Segers, J. (2003) 'Inference for Clusters of Extreme Values', {\it J. R. Statist. Soc. B.}, Vol. 65, pp. 545--556.

\bibitem[\protect\citeauthoryear{Jessen and Mikosch}{2006}]{jessen}
Jessen, A.H. and Mikosch, T. (2006) ' Regularly varying functions', {\it PUBLICATIONS DE L'INSTITUT MATH\'{E}MATIQUE Nouvelle s\'{e}rie}, Vol. 80, No. 94, pp. 171-–192.

\bibitem[\protect\citeauthoryear{L\"{a}ttil\"{a} and Hilmola}{2012}]{Lattila}
L\"{a}ttil\"{a}, L. and Hilmola, O.-P. (2012) 'Forecasting long-term demand of largest Finnish sea ports', {\it Int. J. of Applied Management Science}, Vol. 4, No. 1, pp. 52-–89.

\bibitem[\protect\citeauthoryear{Leadbetter et al.}{1983}]{leadbetter83}
  Leadbetter, M.R.,  Lingren, G. and  Rootz\'{e}n, H. (1983) 'Extremes and Related Properties of Random Sequence and Processes',   {\it Springer, New York}.

\bibitem[\protect\citeauthoryear{Leskovec and Krevl}{2014}]{snapnets}
Leskovec, J. and Krevl, A. (2014) 'SNAP Datasets: Stanford Large Network Dataset Collection'


\bibitem[\protect\citeauthoryear{Linyuan L\"{u} et al.}{2012}]{lu}
Linyuan L\"{u}, Mat\'{u}\v{s} Medo, Chi Ho Yeung, Yi-Cheng Zhang, Zi-Ke Zhang and Tao Zhou (2012)
'Recommender systems', {\it Physics Reports}, Vol. 519, No. 1, pp. 1--49.

\bibitem[\protect\citeauthoryear{Mahajan et al.}{2016}]{Mahajan}
Mahajan, V., Misra, R., Mahajan, R. (2016) 'Review on factors affecting customer churn in telecom sector', {\it Int. J. of Data Analysis Techniques and Strategies}, In Press.

\bibitem[\protect\citeauthoryear{Markovich and Krieger}{2010}]{markovichkrieger}
Markovich N.M., Krieger, U.R. (2010) 'Statistical Analysis and Modeling of Skype VoIP Flows', {\it Computer Communications}, Vol.33, pp. S11--S21.

\bibitem[\protect\citeauthoryear{Markovich}{2013}]{markovich13}
 Markovich, N.M. (2013) 'Quality Assessment of the  Packet Transport  of Peer-to-Peer Video Traffic in High-Speed Networks', {\it Performance Evaluation}, Vol. 70, pp. 28--44.

\bibitem[\protect\citeauthoryear{Markovich}{2014}]{markovich14}
 Markovich, N.M. (2014) 'Modeling clusters of extreme values', {\it Extremes}, Vol. 17, No. 1, pp. 97--125.

\bibitem[\protect\citeauthoryear{Markovich}{2015}]{markovich15}
 Markovich, N.M. (2015) 'Extremes control of complex systems with applications to social networks', {\it 15th IFAC/IEEE/IFIP/IFORS Symposium on Information Control Problems in Manufacturing INCOM 2015},
May 11-13, 2015, Vol. 48, No. 3, Edited by Alexandre Dolgui, Jurek Sasiadek and Marek
Zaremba, pp. 1296--1301, ISSN 2405-8969, \tc{Ottawa, Canada.}

\bibitem[\protect\citeauthoryear{Markovich}{2016a}]{markovich16a}
 Markovich, N.M. (2016a) 'Erratum to: Modeling clusters of extreme values', {\it Extremes}, pp. 1-4, DOI 10.1007/s10687-015-0237-x.

\bibitem[\protect\citeauthoryear{Markovich}{2016b}]{markovich16b}
 Markovich, N.M. (2016b) 'Clusters of extremes: modeling and inferences' (submitted).


\bibitem[\protect\citeauthoryear{O'Brien}{1987}]{obrien}
O'Brien, G.L. (1987) 'Extreme values for stationary and Markov sequences', {\it Ann. Probab.} Vol. 15, No. 1, pp. 281--291.


\bibitem[\protect\citeauthoryear{Robert}{2009}]{robert09}
 Robert, C.Y. (2009) 'Inference for the limiting cluster size distribution of extreme values', {\it The Annals of Statistics}
 Vol. 37, No. 1, pp. 271-–310.

\bibitem[\protect\citeauthoryear{Robert}{2013}]{robert13}
 Robert, C.Y. (2013) 'Automatic declustering of rare events', {\it Biometrika}, Vol. 100, pp. 587--606.

\bibitem[\protect\citeauthoryear{Roberts et al.}{2006}]{roberts}
 Roberts, G. O.,  Rosenthal, J. S.,  Segers, J. and  Sousa B. (2006) 'Extremal Indices, Geometric Ergodicity of  Markov Chains and MCMC',
	{\it Scandinavian Journal of Statistics}, Vol. 9,  pp. 213--229.

\bibitem[\protect\citeauthoryear{Robinson et al.}{2000}]{robinson}
Robinson, M.E., Tawn, J.A. (2000) 'Extremal analysis of processes sampled at different frequences', {\it Journal of the Royal Statistical Society Series B.}, Vol. 62, No. 1, pp. 117-–135.

\bibitem[\protect\citeauthoryear{Smith and Weissman}{1994}]{SmithWeissman}
Smith, R. L. and WEISSMAN, I. (1994) 'Estimating the extremal index', {\it Journal of the Royal Statistical Society Series B.}, Vol. 56, No. 3, pp. 515--528.

\bibitem[\protect\citeauthoryear{Weissman and Novak}{1998}]{WeissmanNovak}
Weissman, I. and Novak, S.Yu. (1998) 'On blocks and runs estimators of the extremal index',  {\it Journal of Statistical Planning and and Inference}, Vol. 66,  pp. 281-–288.

\bibitem[\protect\citeauthoryear{Zhao and Gilbert}{2015}]{Zhao}
Xiaofeng Zhao and Gilbert, K. (2015) 'A statistical control chart for monitoring customer waiting time', {\it Int. J. of Data Analysis Techniques and Strategies}, Vol. 7, No. 3, pp. 301–-321.
\end{thebibliography}
\end{document}